\documentclass{article}

\usepackage{color}
\usepackage{hyperref}  %%[hidelinks] 
\usepackage{graphicx}
\usepackage{amsfonts}
\usepackage{bm}

\makeindex

 2%
\def\Ba{{\bm{\alpha}}}

\def\Bn{{\bm{\nu}}}
\def\Bo{{\bm{\omega}}}
\def\BD{{\bm{\D}}}
\def\BX{{\bm{\X}}}
\def\BDpr{{\bm{\partial}}}

\mathchardef\Bd   = "050E  %delta
\mathchardef\Be   = "0522  %varepsilon
\mathchardef\Bee  = "050F  %epsilon
\mathchardef\Bh   = "0511  %eta
\mathchardef\Bthh = "0512  %teta
\mathchardef\Bth  = "0523  %varteta
\mathchardef\Bi   = "0513  %iota
\mathchardef\Bk   = "0514  %kappa
\mathchardef\Bl   = "0515  %lambda
\mathchardef\Bm   = "0516  %mu
\mathchardef\Bx   = "0518  %xi
\mathchardef\Bom  = "0530  %omi
\mathchardef\Bp   = "0519  %pi
\mathchardef\Br   = "0525  %ro
\mathchardef\Bro  = "051A  %varrho
\mathchardef\Bs   = "051B  %sigma
\mathchardef\Bsi  = "0526  %varsigma
\mathchardef\Bt   = "051C  %tau
\mathchardef\Bu   = "051D  %upsilon
\mathchardef\Bf   = "0527  %phi
\mathchardef\Bch  = "051F  %chi
\mathchardef\Bps  = "0520  %psi
\mathchardef\Bome = "0524  %varomega
\mathchardef\BG   = "0500  %Gamma
\mathchardef\BTh  = "0502  %Theta
\mathchardef\BL   = "0503  %Lambda
\mathchardef\BP   = "0505  %Pi
\mathchardef\BS   = "0506  %Sigma
\mathchardef\BU   = "0507  %Upsilon
\mathchardef\BF   = "0508  %Fi
\mathchardef\BPs  = "0509  %Psi
\mathchardef\BO   = "050A  %Omega
\mathchardef\Bstl = "053F  %*

\newdimen\xshift \newdimen\xwidth \newdimen\yshift \newdimen\ywidth
\def\fline{\hbox to\hsize}
\def\ins#1#2#3{\vbox to0pt{\kern-#2pt\hbox{\kern#1pt #3}\vss}\nointerlineskip}

\def\eqfig#1#2#3#4#5{
\par\xwidth=#1pt \xshift=\hsize \advance\xshift
by-\xwidth \divide\xshift by 2
\yshift=#2pt \divide\yshift by 2
\fline{\hglue\xshift \vbox to #2pt{\vfil
#3 \includegraphics{#4.eps}
}\hfill\raise\yshift\hbox{#5}}}

\def\8{\write12}

\let\a=\alpha    \let\g=\gamma  \let\d=\delta     \let\e=\varepsilon
  \let\h=\eta    \let\th=\theta \let\k=\kappa     \let\l=\lambda
\let\m=\mu    \let\n=\nu          \let\p=\pi        \let\r=\varrho
      
     \let\ch=\chi
\let\G=\Gamma \let\D=\Delta       \let\X=\Xi
     \let\F=\Phi

\let\dpr=\partial
\def\tende#1{\,\vtop{\ialign{##\crcr\rightarrowfill\crcr
 \noalign{\kern-1pt\nointerlineskip} \hskip3.pt${\scriptstyle
 #1}$\hskip3.pt\crcr}}\,}
\def\otto{\,{\kern-1.truept\leftarrow\kern-5.truept\to\kern-1.truept}\,}
\def\fra#1#2{{#1\over#2}}

\def\CC{{\cal C}}\def\DD{{\cal D}}

\def\KK{{\cal K}}

\def\RR{{\cal R}}\def\TT{{\cal T}}

\def\ZZ{{\cal Z}}%
\def\V#1{{\bf#1}}%
\def\defi{\,{\buildrel def\over=}\,}%
\def\*{{\vskip2mm}}\def\0{\noindent}%
\def\lis#1{\overline{#1}}%
\def\iniz{\setcounter{equation}{0}{%
\rhead{\thepage}\lhead{{{{\small\bf\thesection:}\ \SEC\ \  \tiny\today}}}}
\renewcommand{\theequation}{\arabic{section}.\arabic{equation}}%
}
\def\inizA{\setcounter{equation}{0}{%
\rhead{\thepage}\lhead{{{{\small\bf\thesection:}\ \SEC\ \  \tiny\today}}}}
\renewcommand{\theequation}{\Alph{section}.\arabic{equation}}
}

\def\be{\begin{equation}}%
\def\ee{\end{equation}}%

\usepackage{fancyhdr}

\def\wt{\widetilde}\def\tto{\Rightarrow}

\def\eg{{\it e.g.\ }}
\def\ie{{\it i.e.\ }}
\def\tto{\Rightarrow}

\def\alert#1{{\color{ired}#1}}
\definecolor{iblue}{RGB}{65,105,225}
\definecolor{ired}{RGB}{220,20,60}
\definecolor{igreen}{RGB}{50,205,50}
\definecolor{ipurple}{RGB}{75,0,130}
\definecolor{iochre}{RGB}{218,165,32}
\definecolor{iteal}{RGB}{51,204,204}
\definecolor{imauve}{RGB}{204,51,153}

\def\eqalign#1{\null\,\vcenter{\openup\jot
  \ialign{\strut\hfil$\displaystyle{##}$&$\displaystyle{{}##}$\hfil
    \crcr#1\crcr}}\,}

\def\Eq#1{{\label{#1}}}%
\def\equ#1{(\ref{#1})}

\begin{document}

\title{\alert{\bf Quasi periodic Hamiltonian Motions, Scale Invariance,
    Harmonic Oscillators}}

\author{Giovanni Gallavotti
  \\ {\it Dipartimento di Fisica,}
  \\ {\it Universit\`a di Roma-La Sapienza and  I.N.F.N.-Roma1}
  }
\date{}
\maketitle

\begin{abstract}
The work of Kolmogorov, Arnold and Moser appeared just before the
renormalization group approach to statistical mechanics was proposed by
\cite{Wi971a}: it can be classified as a multiscale approach which also
appeared in works on the convergence of Fourier's series,
\cite{Ca966,Fe973}, or construction of Euclidean quantum fields,
\cite{Ne966}, or the scaling analysis of the short scale behaviour of
Navier-Stokes fluids, \cite{CKN982}, to name a few which originated a great
variety of further problems.  In this review the KAM theorem proof will be
presented as a classical renormalization problem with the harmonic
oscillator as a ``trivial'' fixed point.
\end{abstract}

\def\SEC{Introduction}
\section{\SEC}
\label{sec1}\iniz
%\lhead{\SEC}

The KAM theorem can be regarded as a multiscale analysis of the stability
of the harmonic oscillator viewed as a fixed point of a transformation
which enlarges a region of phase space focused around a nonresonant quasi
periodic motion.
The problem considers a Hamiltonian
\be H_0(\V A,\Ba)=\frac12 (\V A\cdot J_0\V A)+\Bo_0\cdot{\V A}
+f_0(\V A,\Ba)\equiv h_0+f_0\Eq{e1.1}\ee
real analytic for $(\V A,\Ba)\in (\DD_\r\times \TT^\ell)$ with:
$\DD_\r=\{\V A\in R^\ell, |A_j|<\r\}$, $\TT^\ell$ the $\ell$-dimensional
torus $[0,2\p]^\ell$ identified with unit circle $\{\V z |\, z_j=e^{i\a_j},
j=1,\ldots,\ell\}$, $\Bo_0\in R^\ell$ and $J_0$ could be a $\ell\times\ell$
{\it non degenerate} symmetric matrix ($\det J_0\ne0$) but here it will be
just the identity matrix time a constant, to simplify notations.

The Hamiltonian is supposed holomorphic in the complex region
$\CC_{\r_0,\k_0}$ with size of the perturbation $f_0$ measured by $\e_0$:
\be\eqalign{
  \CC_{\r_0,\k_0}&\defi\{(\V A,\V
  z)||A_j|\,\le\,\r_0,\ e^{-\k_0}\le| e^{i\a_j}|\le
  e^{\k_0}, j=1,\ldots,\ell\}\subset \CC^{2\ell}\cr
  \e_0&=||\BDpr_{\V A}f_0||_{\r_0,\k_0}
  + \frac1{\r_0} ||\BDpr_{\Ba}f_0||_{\r_0,\k_0},\qquad{\rm with:}\cr
  ||f&||_{\r_0,\k_0}\defi\max_{\CC_{\r_0,\k_0}} |f(\V A,\V z)|,\qquad
  \forall f\ {\rm holomorphic\ in}\ \CC_{\r_0,\k_0} \cr}
\Eq{e1.2}\ee
with $\r_0>0,\k_0>0$, $z_j\equiv e^{i\a_j}$;
generally $\CC_{\r,\k}(\V {\lis A})$ will denote a
polydisk centered at $\V {\lis A}$, \ie defined as in Eq.\equ{e1.2} with
$|A_j-\lis A_j|\le \r$ replacing  $|A_j|\le \r$ and $e^{-\k}\le |z_j|\le
e^\k$; polydisks centered at the ``origin'' will be simply denoted
$\CC_{\r,\k}$ and called ``centered polydisks''.

It is supposed, no loss of generality, that the $\Ba$-average $\lis f_0(\V
A)$ of $f_0(\V A,\Ba)$ vanishes at $\V A=\V0$.

Set $|\V A|=\max |A_j|, |\V z|=\max|z_j|,\ \forall \V A,\V z\in \CC^\ell$.

The idea is to focus attention on the center of $\CC_{\r_0,\k_0}$ where, if
$\e_0=0$, a motion (``free motion'') takes place which is quasi periodic
``with spectrum'' $\Bo_0$. This is done by changing variables in a small
polydisk $C_{\wt\r,\wt\k}(\V a)\subset C_{\r_0,\k_0}$, eccentric if $\V
a\ne\V 0$, that is then recentered and enlarged back to the original size
so that it contains $\CC_{\r_0,\k'_0}$ with $\k'_0>\frac12\k_0$.

The motions developing in the initial polydisk can be studied as ``through
a microscope'': in the good cases (\ie under suitable assumption on the
initial parameters $J_0,\Bo_0$ and $f_0$) the Hamiltonian will turn out to
be substantially closer to that of a harmonic oscillator (described by its
``normal'' Hamiltonian $\Bo_0\cdot \V A$ in the variables $\V A,\Ba)$.

Iterating the process the Hamiltonian changes but, {\it remaining analytic
  in the {\it same polydisk} $C_{\r_0,\frac12\k_0}$}, converges to that of
  a harmonic oscillator: the interpretation will be that, looking very
  carefully in the vicinity of the torus $\TT_{\Bo_0}=\{\V A=\V
  0,\Ba\in[0,2\p]^\ell\}$, also the perturbed Hamiltonian exhibits a
  harmonic motion with spectrum $\Bo_0$: the result, proved below, is the
  KAM theorem.

This is not only reminiscent of the methods called ``renormalization
group'', RG, in quantum field theory but in this review it will be shown to
be just a realization of them, %correcting an error in \cite{Ga986} and
adapting
%the correction
\cite{Ga986} to more recent views on the RG.
%: see concluding remarks where the mentioned error is recalled.

\def\SEC{A formal coordinate change}
\section{\SEC}
\label{sec2}
\iniz

The Hamiltonian Eq.\equ{e1.1}, considered as a holomorphic function on a
domain $\CC_{\r_0,\k_0}$ (Eq.\equ{e1.2}), will be denoted $H_0=h_0+f_0$.
The label $0$ is attached since the beginning because $H_n,f_n,\r_n,\k_n$
will arise later with $n=1,2,\ldots$.

The frequency spectrum $\Bo_0$ will be supposed ``Diophantine'', \ie for
some $C_0>0$ it is, for all
$\V0\ne\Bn\in \ZZ^\ell$ where $\ZZ^\ell$ is the lattice of the integers:

\be |\Bo_0\cdot\Bn|^{-1}\, < \,C_0\, |\Bn|^{\ell},
\quad \forall \Bn\ne\V0\Eq{e2.1}\ee
and the latter inequality will be repeatedly used to define canonical
transformations with generating functions of the form $\F(\V
A',\Ba)+(\V A'+\V a)\cdot\Ba$:

\be \V A=\V A'+\V a+\BDpr_\Ba \F(\V A',\Ba), \quad \Ba'=\Ba+\BDpr_{\V A'}\F(\V
A',\Ba)\Eq{e2.2}\ee 
with the function $\F$ chosen so that in
the new coordinates $(\V A',\Ba')$ the perturbation is {\it weaker}, at the
price that the new coordinates will cover a (much) smaller domain, inside the
$\DD_\r\times \TT^\ell$.  \*

{\it To simplify the notations the functions of $\Ba$ will always be
  implicitly regarded as functions of $z_j=e^{i\a_j}$ whenever referring to
  their holomorphy properties, and without further comments their arguments
  will be written as $\V z$ or $\Ba$, as convenient.}  \*

At first the natural choice for $\F$, {\it temporarily forgetting} the
determination of the domain of definition of the transformation
would be

\be \eqalign{
  \F(\V A',\Ba)&=-\sum_{\V 0\ne\Bn\in \ZZ^\ell}
  \frac{f_{0,\Bn}(\V A'+\V a)}{
    i(\Bo_0\cdot\Bn+((\V A'+\V a)\cdot J_0 \Bn)} e^{i\Bn\cdot\Ba}\cr
\V a&=-J_0^{-1}\BDpr_{\V a}\lis f_0(\V a)\cr}
  \Eq{e2.3}\ee
where  $f_{0,\Bn}(\V A)$
is Fourier's transform of $f_0(\V A,\Ba)$, and $\lis
f_0(\V A')$ denotes the average of $f_0(\V A',\Ba)$ over $\Ba$.
%The symbol
%$\BDpr$ with no further labels means gradient with respect to the argument
%of the function to which it is applied.

Then inserting Eq.\equ{e2.2} into $H_0$ the Hamiltonian is transformed,
setting $J_1=J_0+\frac12\BDpr^2_{\bf a}f_0(\bf a)$,
into

\be\eqalign{(0):\quad&H'(\V A',\Ba')=\frac12 (\V A'\cdot J_1\V
  A')+\Bo_0\cdot\V A'\cr
  (1):\quad&+(\Bo_0+J_0(\V A'+\V a))\cdot \BDpr_{\Ba}\F
    +f_0(\V A'+\V a,\Ba)-\lis f_{\V 0}(\V A'+\V a)\cr
  (2):\quad&+\lis f_0(\V A'+\V a)-\lis f_0(\V a)-\BDpr_{\V a}\lis f_0(\V
    a)\cdot\V A'-\frac12\BDpr^2_{\V a}f_0(\V a)\V A' \V A'\cr
  (3):\quad&+ (f_0(\V A'+\V a+\BDpr_\Ba\F,\Ba)-f_0(\V A'+\V
    a,\Ba))+\frac12\BDpr_\Ba\F\cdot J_0\BDpr_\Ba\F\cr
    }\Eq{e2.4}\ee
where the second of Eq.\equ{e2.3} has been used
and a few terms have been added or subtracted (including free addition or
substraction of constants)  so that:
\begin{itemize}
\item{(0)} The unperturbed Hamiltonian,
\item{(1)} This 
  term vanishes if $\F$ is defined via Eq.\equ{e2.3};
\item{(2)} The term is of $O(\e_0(\V A')^3)$,
  hence it is a higher order term if $|\V A'|$
  is small enough.
\item{(3)} The two terms 
  are {\it formally} of higher order in the size
  $\e_0$ of $f_0$.
\end{itemize}

In a domain in which the transformation Eq.\equ{e2.2} could be defined, the
motions would be described by a simpler Hamiltonian which is still the same
integrable Hamiltonian plus a perturbation of higher order in $\e_0$.

{\it However} to make sense of the transformation in Eq.\equ{e2.2} it is
not only necessary to restrict the variables $(\V A',\Ba)$ to a smaller
domain, but it has to be possible to solve the implicit functions problem
in Eq.\equ{e2.2},\equ{e2.3} (namely to express $(\V A,\Ba)$ in terms of
$(\V A',\Ba')$ and viceversa, and finding $\V a$), but also the denominator
in Eq.\equ{e2.3} will have to be modified to avoid dividing by $0$: which
will happen, for generic $f_0$ and for some $\Bn$, on a dense set of $\V
A'\in \DD_{\r_0}$, if $J_0$ is not singular (as it is being supposed).
Therefore the map in Eq.\equ{e2.2} will now be modified and defined
properly after recalling the notion of dimensional estimate.

\def\SEC{Dimensional estimates}
\section{\SEC}
\label{sec3}
\iniz

The very nature of the stability of quasi periodic motions is that it is a
multiscale problem: like many other problems in analysis, from the almost
everywhere convergence of Fourier series of $L_2([0,2\p])$-functions
(\cite{Fe973}), to the study of the possible singularities of the
Navier-Stokes problem (\cite{CKN982}), to the convergence of the functional
integrals arising in quantum field theory (\cite{Ga985b}), to name a
few. The {\it renormalization group} method, \cite{WK974,BG995}, unifies
the approaches developed to study such problems.

The main feature of the renormalization group applications is  their
being based on what will be called here ``{\it dimensional estimates}''.

Dimensional estimates deal with elementary bounds on holomorphic functions.
Let $g(z)$ be any holomorphic function in a closed domain $C\subset\CC$
(domain $\tto$ closure of an open set in the complex plane $\CC$). The
function $g$ can be bounded, toghether with its Taylor coefficients, in
terms of $||g||_C=\max_{z\in C}|g(z)|$, inside the region $C_\d$ consisting
of the points in $C$ at distance $\ge \d$ from the boundary of $C$:

\be |\dpr^n_z g(z)|\le\,  n!\, ||g||_{\lower 1mm \hbox{$\scriptstyle C$}}\,
\d^{-n},\qquad \forall z\in
\CC_\d,\ n\ge0\Eq{e3.1}\ee
A consequence is that if $g$ is holomorphic in a disk
$C_\r=\{z|\,|z|\le\r\}$ or in an poli-annulus $\G_\k=\{\V z| e^{-\k}\le
|z_j|\le e^\k,\, j=1,\ldots,\ell\}$ then the following elementary bounds on
the derivatives of $g$ or, respectively, the Fourier coefficients $g_\Bn$ of
the function $g(e^{i\a})$ hold

\be \eqalign{
  &||\dpr^n_z g||_{\lower 1mm \hbox{$\scriptstyle C_{\r'}$}}\le
  \,  n!\, ||g||_{\lower 1mm \hbox{$\scriptstyle C_{\r}$}}\,  (\r-\r')^{-n},
  \qquad \forall n\ge0
  \cr
  & |g_{\Bn}|\le\,
  ||g||_{\lower 1mm \hbox{$\scriptstyle \G_{\k}$}}\,  e^{-\k |\Bn|},
  \qquad\forall \Bn\in\ZZ,\ |\Bn|=\sum_{i=1}^\ell |\n_i|\cr}
\Eq{e3.2}
\ee

Holomorphic functions $g$  of $\ell$ or $2\ell$ arguments
will be considered, in the following, in domains

\be \eqalign{ &\CC_\r=\{\V A | |A_j|\le\r, j=1,..,\ell\},\quad
  \G_\k=\{\V z| e^{-\k}\le |z_j|\le e^\k, j=1,..,\ell\}\cr
  &\CC_{\r,\k}=\CC_\r\times\G_\k\cr } \Eq{e3.3}\ee
and their maxima will be denoted by appending labels $\r$ or $\k$ or
$\r,\k$, as appropriate, to the symbol $||g||$.

Hence if
$||g||_{\r,\k}=\e$ the bounds

\be \eqalign{
  &||g_{\Bn}||_\r \le \,\e\, e^{-\k|\Bn|},\kern2cm  \forall
  \Bn\in\ZZ^\ell, \V A \in \CC_{\r'}\cr& 
  ||\BDpr_{\V A}^n g_{\Bn}||_{\r',\k}\le\,  n!\, \e \, e^{-\k|\Bn|}
  (\r-\r')^{-n},\quad  \forall
  \Bn\in\ZZ^\ell, \V A \in \CC_{\r'}\cr
  \cr}
\Eq{e3.4}
\ee
hold and will be called {\it dimensional bounds}.

Summarizing: the dimensional bounds say that the $n$-th derivatives of a
function holomorphic in a domain $C$ are bounded, at a point $z$ at distance
$\d$ from the boundary of $C$, by the maximum of the function in $C$ divided
by the $n$-th power of the distance of $z$ to the boundary $\dpr C$ of $C$
times $n!$ (``Cauchy's theorem'').

In the following essentially all bounds will be ``dimensional'': and each
new bound presented may contain some new constants labeled
$c_i,\g_i$; such constants will only depend on the number of degrees of
freedom $\ell$ and, for simplicity, will be chosen so that $c_i\le
c_{i+1},\g_i\le\g_{i+1}$ and .

\def\SEC{A canonical map}
\section{\SEC}
\label{sec4}
\iniz

The ``renormalization group'' is
a map $\RR$ whose iterations can be
interpreted as successive magnifications zooming on ever smaller regions of
phase space in which motions develop closer and closer to the searched
quasi periodic motion of spectrum $\Bo_0$. 

At step $n=0,1,\ldots$ the motions will be described by a
Hamiltonian $H_n+f_n$ which will be the sum of three terms

\be \frac12{\V A\cdot
  J_{n}\V A}+\Bo_0\cdot \V A+f_n(\V A,\Ba),\Eq{e4.1}\ee
see Eq.\equ{e1.1}. In the renormalization group nomenclature and {\it under
  the conditions Eq.\equ{e2.1} and $\det J_0\ne0$} the first and third
terms would be called
``{\it
  irrelevant}'' and the intermediate (\ie the normal form for the
$\ell$-dimensional harmonic oscillators Hamiltonian) would be called a
``{\it marginal
  trivial fixed point}'': the reason behind the latter names will be
be mentioned below.

Introducing the parameters $\e_n,J_n,\r_n,C_n,\k_n$, characterizing $H_n$
in the same sense in which $\e_0,J_0,\r_0,C_0,\k_0$ characterize $H_0$ in
Eq.\equ{e1.2}, it is convenient, for the purpose of a rapid evalutation of
several estimates, to keep in mind that the following ``dimensionless''
quantities,

\be \h_n=\e_nC_n,\qquad \ e^{\k_n}\Eq{e4.2}\ee
will naturally occurr in the dimensional estimates: the latter will,
therefore, be expressed as products of selected dimensionless quantities
times a suitable factor chosen among the dimensional parameters
$\e_n,\r_n,C_n,J_n$.  

{\it All bounds} will be carefully written so that they will involve only
  dimensionless constants and, when needed, a factor to fix the dimensions.
  Furthermore the construction of the sequence $H_n$ will be so designed
  that
  \be C_n\equiv C_0,\r_n\ll\r_0,\k_n=\k_{n-1}-4\d_n> \frac12\k_0\Eq{e4.3}\ee
  with $\d_n$
  defined so that $\k_0\ge\k_n\ge\frac12\k_0$; to fix the ideas $\d_n$ will be
  fixed as $\d_n= (n+10)^{-2}\k_0$, the size $\e_n$ of  $f_n$ will tend to $0$
  provided $\e_0$ is small enough, while $J_n=J_0$.

It will not be restrictive to suppose, {\it initially}:

\be 
C_0\r_0 J_0<1,\qquad 2^{-1}<e^{\frac{\k_0}2}<e^{\k_n}<\,e^{\k_0}<2 \Eq{e4.4}\ee
because the theorem will apply for $\e_0$ small enough and $\r_0,\k_0$ can
be {\it initially} restricted as needed. Furthermore it is important to
keep in mind that the bounds that follow are naive dimensional bounds
derived {\it without any optimization attempt}, yet they will suffice for a
complete proof.

To define properly a tranformation inspired by Eq.\equ{e2.2} and to
eliminate the mentioned possible divisions by $0$, {\it while still keeping
  $H'$ in Eq.\equ{e2.4} formally close to $H_0$} as in Sec.\ref{sec2}, the
first task is to determine the 
shift $\V a$, Eq.\equ{e2.3}.

The implicit equation Eq.\equ{e2.3} for $\V a$, $\V a=-J_0^{-1}\BDpr_{\V a} \lis
f_0(\V a)\defi \V a_0+\V n(\V a)$, with $\V a_0=-J_0^{-1}\BDpr_{\V a}
\lis f_0(\V
0)$ can be solved under a smallness condition on $\e_0$ obtaining $\V a$
close to $\V a_0$.

This follows from an application of a general implicit function theorem
yielding the existence of a constant $\ch$ such that the smallness
condition $|\V n|_{\r_0}\r_0^{-1}<\ch$ implies existence of a
solution. Since $|\V n|_{\r_0}\r_0^{-1}$ is dimensionally bounded by
$\e_0J_0^{-1}\r_0^{-1}\defi\th_0$, a condition for the solubility of the
equation is:
\be\th_0=\e_0\r_0^{-1}J_0^{-1}<\ch\qquad\tto\qquad |\V
a|<\th_0\r_0<\ch\r_0<\frac1{16}\r_0
\Eq{e4.5}\ee
The choice $\ch=\frac1{16}$ (implied in general by the estimate of $\ch$
reproduced for completeness in  Appendix \ref{secA} below) is useful
for the coming analysis (with no attention to an optimal $\ch$-value).

The function $f_0(\V A'+\V a,\Ba)$ will then be defined and analytic
in $C_{\frac34\r_0,\k_0}$ (from $\frac34+\frac1{16}<1$). 
Then proceed to build $\F$, but replace Eq.\equ{e2.3} with its second
order expansion in $J_0$:

\be\eqalign{\F_0(\V A',\Ba)=&-\sum_{\V 0\ne\Bn\in \ZZ^\ell}
  \frac{f_{0,\Bn}(\V A'+\V a) e^{i\Ba\cdot\Bn}}{i\Bo_0\cdot\Bn} \cr&
  \cdot\Big(1-\frac{J_0(\V A'+\V a)\cdot\Bn}{\Bo_0\cdot\Bn}+ (\frac{J_0(\V
    A'+\V a)\cdot\Bn}{\Bo_0\cdot\Bn})^2\Big) \cr}\Eq{e4.6}\ee
The function $\F_0$ is well defined in the
polydisk $\CC_{\frac34\r_0,\k_0-\d_0}$ as seen via the following general
dimensional bounds (given in Eq.\equ{e3.4} on functions bounded by $\e_0$
and holomorphic in a domain $\CC_{\r_0,\k_0}$).

Taking into account the Diophantine inequality Eq.\equ{e2.1}, for $0\le
\d_0<\k_0$, the definitions Eq.\equ{e4.2},\equ{e4.3},\equ{e4.4} and the
dimensional inequality Eq.\equ{e3.4}, with the restrictions Eq.\equ{e4.4},
leads to:
\be \eqalign{
  ||\F_0||_{\frac34\r_0,\k_0-\d_0}\le& \e_0\r_0 \sum_{\Bn\ne\V0}
\frac{e^{-\d_0|\Bn|}}{|\Bo_0\cdot\Bn|}(1+
\frac{|J_0|\r_0|\Bn|}{|\Bo_0\cdot\Bn|}\cr
&+(\frac{|J_0|\r_0|\Bn|}{|\Bo_0\cdot\Bn|})^2)
\le \g_1 \h_0\r_0 \,\d_0^{-c_1}
\cr
|\F_{0,\Bn}(\V A')|\le& \g_1 \h_0\,\r_0 \,\d_0^{-c_1}\,
e^{-\k_0|\Bn|},\qquad \forall |\V
A'|<\r_0
\cr}
\Eq{e4.7}\ee
with $\g_1,c_1$ are dimensionless constant (depending only on the number of
degrees of freedom $\ell$, \eg $c_1=5\ell+2$), $\h_0=\e_0 C_0$ and
$J_0C_0\r_0<1$ have been used.

Hence the functions in the {\it r.h.s} of Eq.\equ{e2.2} admit the
dimensional bounds:
\be\kern-4mm\eqalign{
 & ||\BDpr_{\Ba}\F_0||_{\frac23\r_0,\k_0-2\d_0}
  \le \g_2 \h_0\,\r_0 \d_0^{-c_2},\ 
%  \cr &
  ||\BDpr_{\V A'}\F_0||_{\frac23\r_0,\k_0-2\d_0}
  \le \g_3 \h_0 \,\d_0^{-c_3}
  \cr&
  ||\BDpr^2_{ \Ba \V A'}\F_0||_{\frac23\r_0,\k_0-2\d_0}  
  \le \g_4 \h_0\,\d_0^{-c_4},\ 
%  \cr&
  ||\BDpr^2_{\V A' \V A'}\F_0||_{\frac23\r_0,\k_0-2\d_0}
\le \g_5 \h_0 \r_0^{-1}\d_0^{-c_5}\cr
}
\Eq{e4.8}\ee
where the derivatives with respect to $\a_j$ should be interpreted as $i
z_j\dpr_{z_j}$ for $z_j=e^{i\a_j}$ in the domain $(\V A',\Ba)\in
C_{\frac23\r_0,\k_0-2\d_0}$, and the constants $\g_i,c_i$
can be fixed to depend only on $\ell$.
The radius is reduced to $\frac23\r_0$ to allow simple dimensional bounds
using $\frac34-\frac1{16}>\frac23$ (taking into account the second
inequality in Eq.\equ{e4.5}).

To define the canonical transformation $(\V A',\Ba')\to(\V A,\Ba)$ the
implicit functions in Eq.\equ{e2.2} have to be solved. This can be done
quite easily if one is willing to define the map only for $(\V A',\Ba')$
contained in a small enough domain.

The condition to express $(\V A,\Ba)$ in terms of $(\V A',\Ba')\in
\CC_{\r',\k'}$ with $\r'=\frac12\r_0,\k'=\k_0-3\d_0$ is prescribed via an
implicit function theorem for analytic functions, see for instance
propositions 20,21 in Sec.5.11 and Appendix N in \cite{Ga983}, or
\cite[Appendix3]{Ga985}.

The theorem is proved following the lines of the
analogous result ``for disks'' leading to Eq.\equ{e4.5} (discussed in
Appendix A below) adapting it to polydisks and the condition is obtained on
dimensional grounds as the bound (on the Jacobian of the implicit equations
Eq.\equ{e2.2})

\be \eqalign{
  &||\BDpr^2_{\V A' \Ba}\F_0||_{\frac23\r_0,\k_0-2\d_0}<
  \g_6 \h_0\,\d_0^{-c_6}<1\cr
  }\Eq{e4.9}\ee
where the first inequality is just the bound Eq.\equ{e4.8} on the {\it
  l.h.s.} with $\g_4$ modified into a larger $\g_6$ and $c_6$ a constant
(\eg $5\ell+4$).
 
This can be obtained, again reducing the radius from $\frac23\r_0$
to $\frac12\r_0$ for ease of dimensional bounds, by first
fixing $\V A'\in C_{\frac12\r_0}$ so that the second inequality in
Eq.\equ{e4.9} simply {\it implies injectivity} of the map
$\Ba'=\Ba+\dpr_{\V A'}\F_0(\V A',\Ba)$ for $\Ba\in C_{\k_0-2\d_0}$, for all
$\V A'$ fixed in $C_{\frac12\r_0}$; it implies also $\Ba\in C_{\k_0-2\d_0}$
for $\Ba'\in C_{\k_0-3\d_0}$ if $\g_6$ is large enough, see appendix
\ref{secB}. Therefore, given $\V A'\in C_{\frac12\r_0}$ and using the
injectivity, $\Ba$ can be computed from $\Ba'$ in the form

\be\eqalign{
  &\Ba=\Ba'+\BD(\V A',\Ba'), \qquad \Ba'\in C_{\k_0-3\d_0}, \forall \V
  A'\in C_{\frac12\r_0}\cr
  &\BD(\V A',\Ba')\equiv -\dpr_{\V A'}\F_0(\V A',\Ba)\cr
  &||\BD||_{\frac12\r_0,\k_0-3\d_0}<\g_7\h_0\d_0^{-c_7}<\d_0,
  \cr
}
\Eq{e4.10}\ee
where the second line in Eq.\equ{e4.10} is an identity which implies, via
Eqs.\equ{e4.8},\equ{e4.9}, the inequalities in the third line, where
$\g_7,c_7$ are suitable positve constants.

The second inequality in Eq.\equ{e4.9} also insures the injectivity of $\V
A=\V A'+\dpr_\Ba\F_0(\V A',\Ba)$ for $\V A'$ in $C_{\frac12\r_0}$, for all
$\Ba$ fixed in $C_{\k_0-2\d_0}$, therefore for all $\Ba'$ in
$C_{\k_0-3\d_0}$.

Hence $\BD(\V A',\Ba')$ is defined in $C_{\frac12\r_0,\k_0-3\d_0}$ and the
  angles $\Ba$ can be expressed in terms of $\Ba',\V A'$; and it is
  possible to express, for each $\Ba'\in C_{\k_0-3\d_0}$, $\V A$ in terms
  of $\V A',\ \forall \V A'\in C_{\frac12\r_0}$: simply by substituting
  $\Ba$ by $\Ba'+\BD(\V A',\Ba')$ to find:
\be\eqalign{
  &\V A=\V A'+\V a+\BX(\V A',\Ba')\cr
  &\BX(\V A',\Ba')\equiv \dpr_\Ba\F_0\Big(\V A',\Ba'+\BD(\V
  A',\Ba')\Big),\cr
  }
\Eq{e4.11}\ee
For $(\V A',\Ba')\in C_{\frac1{2}\r_0,\k_0-3\d_0}$ the $(\V A,\Ba)$ will
vary inside the original domain.

Then again Eqs.\equ{e4.8},\equ{e4.9}, if $\g_7 \h_0\d_0^{-c_7}<1$
for $\g_7$ suitably larger than $\g_2$ and $c_7=c_2$, yield

\be|\V A|<\frac12 \r_0+||\BX||_{\frac12\r_0,\k_0-3\d_0}\le
  \frac12 \r_0+\g_2 \h_0\r_0 \d_0^{-c_2}<\frac34\r_0\Eq{e4.12}\ee

  Collecting all conditions to define $\V a,\BD,\BX$ a canonical map

\be\eqalign{ \V A=&\V A'+\V a+\BX(\V A',\Ba'),
  \qquad \Ba=\Ba'+\BD(\V A',\Ba')\cr
  ||\BX&||_{\frac12 \r_0,\k_0-3\d_0}< \g_8 \h_0\,\r_0\d_0^{-c_8} \cr
  ||\BD&||_{\frac12 \r_0,\k_0-3\d_0}< \g_8 \h_0\,  \d_0^{-c_8}
  \cr
  } \Eq{e4.13}\ee
will be defined, for suitably chosen $\g_8,c_8$, changing $(\V A',\Ba')\in
C_{\frac12\r_0,\k_0-3\d_0}$ into $(\V A,\Ba)\in C_{\frac34\r_0,\k_0-\d_0}$.

The perturbation function $f_0-\lis f_0(\V a)$ becomes in the new
coordinates $f'_0(\V A',\Ba')=f_0(\V A'+\V a+\BX(\V A',\Ba'),\Ba'+\BD(\V
A',\Ba'))$ and the new Hamiltonian can expressed by replacing $\Ba$ with
$\Ba'+\BD(\V A',\Ba')$ in the three terms in Eq.\equ{e2.4}. This is
discussed in the next section in terms of $\h_0,\d_0$; the conditions
imposed, so far, on the construction can be all implied by the conditions

\be \eqalign{
  &C_0\r_0 J_0<1,\quad
  e^{\k_0}<2\qquad \hbox{initial restrictions}  \cr
  &\e_0 J_0^{-1}\r_0^{-1}<\ch,\qquad \hbox{to define}\qquad
  \V a=-J_0^{-1}\BDpr_{\V a}\lis f_0(\V a) \cr
  &\g_9 \h_0 \d_0^{-c_9}<1,\kern9mm\hbox{to define $\BD,\BX$} \cr}
\Eq{e4.14}\ee
for $\g_9,c_9$ large enough and $\ch$ small enough, see Eq.\equ{e4.5}.

The domain of variability in the initial variables $(\V A,\Ba)$, where the
canonical map is defined, will now contain (at least) a small domain of
shape close to a polydisk ({\it eccentric} because of the 
translation by $\V a$) inside the
initial domain $C_{\r_0,\k_0}$ of the Hamiltonian $H_0$. The small
eccentric polydisk is the image of a {\it centered} polydisk
$\CC_{\frac12\r_0,\k_0-3\d_0}$ in the new variables $(\V A',\Ba')$.

\def\SEC{Renormalization}
\section{\SEC}
\label{sec5}
\iniz

The Hamiltonian $H_0+f_0$ in the new coordinates $\V A',\Ba'$ becomes:

\vskip-5mm
\be\eqalign{
  H'&(\V A',\Ba')=\frac12\V A'\cdot J_1\V A'+\Bo_0\cdot\V A'+f', \qquad
  (\V A',\Ba')\in\CC_{\frac12\r_0,\k_0-3\d_0}
\cr}\Eq{e5.1}\ee

\vskip-3mm \0in the domain $(\V A',\Ba')\in
C_{\frac12\r_0,\k_0-3\d_0}$. The function $f'$ is defined, in the mixed
variables $(\V A',\Ba)$, by Eq.\equ{e2.4}.

\*
\begin{itemize}
\item The contribution 1) in Eq.\equ{e2.4}, {\it
  does not vanish: but it carries the key cancellation} showing that the
sum of  terms individually formally $O(\e_0)$ is in fact of higher
order in $\e_0$ as can be seen via the Fourier's transform of $f_{0}-\lis
f_0=\sum_{\V0\ne \Bn} f_{0,\Bn} e^{i\Ba\cdot\Bn}$ which, after a few
simplifications, is:
\be\eqalign{
  F&\defi(\Bo_0+J_0 (\V A'+\V a)\cdot\BDpr_\Ba\F_0+
  f_0(\V A'+\V a,\Ba)-\lis f_0(\V A'+\V a))\cr
  &=
\sum_{\V0\ne \Bn} f_{0,\Bn}(\V A'+\V a)\frac{(J_0(\V A'+\V a)\cdot\Bn)^3}{(\Bo_0
  \cdot\Bn)^3}e^{i\Ba\cdot\Bn}\cr} \Eq{e5.2}\ee 
\item if $|\V A'|<\wt\r$, $F$ admits a dimensional bound {\it in the sense
  of Eq.\equ{e1.2}}, \ie on $\|\BDpr_{\V A'}
  F\|_{\wt\r,\k_0-3\d_0}+\frac1{\wt\r} \|\BDpr_{\Ba'}
  F\|_{\wt\r,\k_0-3\d_0}$.  Using $J_0\V a=-\BDpr f_0(\V a)$,
  $J_0C_0\r_0<1$, and Eq.\equ{e3.4} together with the bound $|f_{0,\Bn}|<
  \e_0\r_0 e^{-\k_0|\Bn|}$ (derived from $\frac1{\r_0}|\BDpr_{\Ba} f_0|\le
  \e_0$) a dimensional bound on $F$ follows as:
\be \eqalign{
  \le&\ \g_{10}\e_0\,\frac{\r_0}{\wt\r}
    ((J_0\wt\r C_0)^3+(C_0\e_0)^3)\d_0^{-c_{10}}
  \cr&
  =
  \g_{10} 
  \e_0\,\frac{\r_0}{\wt\r}(\frac{\wt\r^3}{\r_0^3}+
  \h_0^3)\d_0^{-c_{10}}<\g_{11}\e_0\h_0^{2(1-\l)}\d_0^{-c_{11}}
  \cr}\Eq{e5.3}\ee
in the polydisk $\CC_{\wt\r ,\k_0-3\d_0}$. The $\wt\r$ will be determined
as $\wt\r=\h_0^{1-\l} \r_0$ with $0<\l<1$ so that 
$\wt\r<\fra12\r_0$ provided $\h_0$ is small enough.

\item The contribution 2) in Eq.\equ{e2.4}, is bounded, {\it still in the
  sense of Eq.\equ{e1.2}}, in a disk of radius $\wt\r=\h_0^{1-\l}\r_0$,
  if, as above, $\h_0$ is small enough, by
\be \g_{13}\e_0\h_0^{2(1-\l)}\Eq{e5.4}\ee
making use of its $\Ba$-independence, which permits to estimate
dimensionally the second
derivative of $\lis f_0(\V A+\V a)$ in a disk of radius $\frac12\r_0$
(rather than of radius $\wt\r_0$): thus it
also yields a contribution to the higher order terms.
\item The terms in the contribution 3) are also dimensionally bounded, {\it
  still in the sense of Eq.\equ{e1.2}}, by:

\be\eqalign{
  &\g_{14}\e_0\h_0^{\l}\d_0^{-c_{14}}\cr}\Eq{e5.5}\ee
in the polydisk $\CC_{\wt\r,\k_0-3\d_0}$, using $C_0\r_0 J_0<1$.
\end{itemize}
\*

Adding the bounds Eq.\equ{e5.3},\equ{e5.4},\equ{e5.5} it is,
for $\l=\fra23$ (\ie $2(1-\l)=\l$):
\be\e_1=(|\BDpr_{\V A'} f'|_{\wt\r,\k_0-4\d_0}+
\wt\r^{-1}|\BDpr_{\Ba'} f'|_{\wt\r,\k_0-4\d_0}<
\g\e_0 \h_0^{2(1-\l)}\d_0^{-c}) \Eq{e5.6}\ee
for $\g,c>0$ suitably fixed, if $C_0\r_0 J_0<1$, (see also Eq.\equ{e4.2}).

A further dimensional estimate on the matrix $\BDpr^2_a \lis f_0(\V a)$ in
Eq.\equ{e2.4} is $J_1<J_0(1+c\th_0),\,J_1>J_0(1-c\th_0)$ possibly increasing the
$c$ appearing in eq.\equ{e5.6}.

The result is that in the coordinates $(\V A',\Ba')\defi(\V A_1,\Ba_1)$ the
motion is Hamiltonian with Hamiltonian $H_1$; and recalling the definitions of
the dimensionless quantities in Eq.\equ{e4.2},\equ{e4.5}:

\be\eqalign{
  &H_1=\frac12\V A_1\cdot J_1\V A_1+\Bo_0\cdot{\V
      A}_1+f_1(\V A_1,\Ba_1)\cr
  & \r_1=\r_0 \h_0^{1-\l},\quad 
  \k_1=\k_0 -\lis\d_0,\quad C_1=C_0\cr
  &\h_1=\g \h_0^{3-2\l}\d_0^{-c},%
  \qquad \th_1= \g\th_0\h_0^{1-\l}\d_0^{-c},
  \qquad J_1 C_1\r_1<J_0 C_0 \r_0<1
  \cr&J_0(1-c\th_0)\,<\,J_1\,<\,J_0(1+c\th_0)\cr}
  \Eq{e5.7}\ee
  where $\g,c$ are constants, $\lis \d_n=4\d_n=\k_0 (n+10)^{-2}$ and
  $\l=\fra23$. 

  The above transformation of coordinates $(\V A,\Ba)\to(\V A_1,\Ba_1)$,
  which will be denoted $\KK_0$, is well defined and holomorphic in the
  domain $C_{\frac12\r_0,\k_0-3\d_0}$ whose $\KK_0$-image contains the
  small polydisk $C_{\wt\r,\k_0-4\d_0}$ {\it provided $\e_0$ is small
    enough} so that the conditions imposed during the construction, namely
  Eq.\equ{e4.14}, and the ones following it, are satisfied and {\it remain
    satisfied under iteration} allowing to define the sequence of maps
  $\KK_n, \,n\ge0$.Because, if $\h_0$ (\ie $\e_0$) is small enough, the map
  in Eq.\equ{e5.7} generates a sequence with $C_n\e_n=C_0\e_n=\h_n$ tending
  to $0$, fixed arbitrarlily $\m\in (0,\frac23)$ and a corresponding
  suitable constant $\lis\g$, superexponentially with

\be \h_n\sim
(\lis\g \h_0)^{(1+\m)^n}, \quad \lis\g>0,\ 0<\m<\frac23,\qquad
J_n<2J_0,\ J_n>\frac12J_0\Eq{e5.8}\ee
and $\th_n$ also tend to $0$ at similar rates (\eg $\th_n\sim$ $(\lis\g'
\h_0)^{(1+\m)^n c'}$), as can be checked by induction
from Eq.\equ{e5.7} with suitable $c'<1,\lis\g'$. This implies that for
all $n\ge0$ the transformations $\KK_n$ can be defined if $\e_0$ (\ie its
dimensionless version $\h_0$) is small enough.

Furthermore $\KK_n$ is seen from Eq.\equ{e4.13} to be close to the identity
within $\g_8 \h_n\d_n^{-c_8}$.  Hence the iteration of the renormalization
procedure defines a sequence of transformations $\KK_n$ under the only
initial condition in Eq.\equ{e4.14} with $\g_9,c_9,\ch^{-1}$ large enough.
 
In the polydisk $\CC_{\r_n,\k_n}$ the motions starting with $\V A_n=\V0$
and (say)$\Ba=\V0$ become closer and closer to the motion of a harmonic
oscillator with frequency spectrum $\Bo_0$ and in the limit $n\to\infty$
all motions in the ``polydisk'' (degenerated to a torus $\V 0\times
\TT^\ell$) are harmonic with spectrum $\Bo_0$. This
is checked simply by remarking that the motion of the initial data is, if
observed in an {\it arbitrarily fixed time $t$}, is superxponentially close to
the harmonic motion $\V A=\V 0,\Ba(t)=\Ba+\Bo_0t$. The torus on which the
motion is quasi periodic is the limit of the tori with equations
$\V A=\V a_n+\BX_n(\V a_n,\Ba'), \Ba=\Ba'+\BD_n(\V a_n,\Ba')$ which is the
torus which at the $n$-th iteration of the renormalization has coordinates
$\V a_n,\Ba')$. The successive corrections to $\V a_n$ and to the functions
$\BX_n,\BD_n$ tend to $0$ superxponentially and their limits
\be \V a_\infty, \ \BX_\infty(\Ba'),\ \BD_\infty(\Ba'), \qquad \Ba'\in
T^\ell
\Eq{e5.9}\ee
define an invariant torus on which motion is $\Ba'\to \Ba'+\Bo_0 t$.

It is also possible to define a sequence of maps $\wt \KK_n$ defined in the
{\it fixed domain} $\CC_{\frac12\r_0,\frac12\k_0}$ by rescaling the
polydisks by a factor $\h_{n-1}^{\frac12}=\r_{n-1}/\r_n$, $n\ge1$ so that
they are all turned into $\CC_{\frac12\r_0,\frac12\k_0}$: the rescaling
transformation will change $\V A_n$ into $\V A'_n=\h_n^{-\frac13}\V A_n$
and the Hamiltonian into

\be\wt H_n=\Bo_0\cdot\V A'_n+ \h_n^{\frac13} \frac13(\V A'_n\cdot J_n\V
  A'_n)
  +\h_n^{-\frac13} f_n(\h_n^{\frac13}\V A'_n,\Ba'_n)
  \tende{n\to\infty} \Bo_0\cdot\V A_{\infty}\Eq{e5.10}\ee
and in the rescaled variables the sizes of the anharmonic terms tend to $0$
superexponentially, taking into account the recursion defined in
Eq.\equ{e5.7} (and that the size of $f_n$ is $\h_n$).

This shows that the perturbation $f_0$ and the twist $J_0$ are, after
renormalization, ``irrelevant operators'' (in Eq.\equ{e5.10} both tend to
$0$ as , while the harmonic oscillator is
a ``fixed point'': in some sense the transformation has the harmonic
oscillator as an {\it attractive fixed point}. This completes a proof of
the KAM theorem, {\it interpreted in the Renormalization Group frame}
\cite{Ko954,Ar963b,Mo962,BGGS984}: it can be  classified as a
``super-renormalizable'' problem, as it requires only a second order
perturbation analysis, Eq.\equ{e4.6}, around the trivial fixed point.  \*

\0{\it Remarks:} (1) { a simpler analysis {\rm (and an instructive warm-up
  exercise)} can be carried also if $J_0=0$ {\rm provided the
  perturbation depends only on the angles} $\Ba$. The independence of
$f_0$ from $\V A$ has the consequence, in the proof development, that all
terms appearing to involve $J_0^{-1}$ actually do not arise at all ({\it but}
the system is elementarily integrable).
\\
(2) The condition $\det J_0\ne0$ is called ``anisochrony condition'' or
``twist condition'': the size $\ne0$ of $\det J_0$ plays an important role
in the above analysis. However invariant diophantine tori, may in certain
cases, exist just for $\e_0$ smaller than a quantity independent on the
size of $\det J_0$; such tori are called ``twistless'', because they can be
shown to exist without invoking the twist condition. This happens in cases
in which $f_0$ depends on $\Ba$ only: and the tori can be constructed via a
simple graphical algorithm, \cite{Ga994b}.
%as first pointed out in \cite[Vol.1,p.133]{Th983}.
The graphical algorithm led, in the twistless cases, also to a new
``direct'' proof of the KAM theorem, \cite{El996,GG995,GBG004}, that was
later extended to the general case, \cite{GM004}.\\
(3) The estimates in the above analysis are far from optimal and optimization
is desirable.}

\def\SEC{Comments}
\section{\SEC}
\label{sec7}
\iniz

The analysis in Sec.\ref{sec5} is a reformulation of the original proof by
Kolmogorov, \cite{Ko954}, reproduced in full detail in \cite{BGGS984} and
used to build a rigorous computation algorithm in \cite{GL999}. The feature
of the approach, common also to Moser's work, \cite{Mo962}, is to use
canonical maps with fixed small denominators: this avoids dealing with $\V
A$ dependent divisors appearing in \cite[p.105]{Ar963b}, reproduced in
\cite{Ga983}.

The renormalization group interpretation has been proposed in in
\cite{Mk983} with prefixed divisors and \cite{Ga982b,Ga986} still dealing
with $\V A$-dependent divisors: the approach developed in Sections 4,5 is
inspired by the latter development but avoids $\V A$-dependent divisors,
hence it is close to \cite{Ko954,Mo962,BGGS984,GL999,CJB999,HI004} and several
other approaches. The definition of $\e_n$, see
Eq.\equ{e1.2}, can be replaced by $\e_0=\max_\CC|f_0|$: this choice would be
possible, jsut with obvious notational changes.\*

The relation between the KAM theorem and the renormalization group has been
used in various forms for its proof, in several papers,
for instance \cite{Ga982b,Mk983,Ga995d,Ga986,BGK999,CGJ998,CJB999,Ko004,Ge010}.

The difference between the approach of Kolmogorov and Moser, with respect
to Arnold's, \cite{Ar963b}, is that in the second the small divisors are
$\V A$-dependent and are controlled by an increasing sequence of cut-offs
on $\Bn$, at each order of the perturbation expansion.

The analysis of the singularity at $\e_0=0$, in the case of resonant quasi
periodic motions (\ie motions which dwell on lower dimensional tori), can
also be pursued via multiscale methods conveniently interpreted as methods
of performing the resummations of the perturbative series, which unlike the
KAM case, are divergent power series, \cite{GG002,GG005e,CGGG006}.
%\end{document}

\renewcommand{\thesection}{\Alph{section}}
\renewcommand{\theequation}{\Alph{section}.\arabic{equation}}

\appendix
\setcounter{equation}{0}
%\begin{appendices}

\def\SEC{Implicit functions in (4.5) (and (4.10),(4.11))}
\section{\SEC}
\label{secA}
\inizA

This appendix is presented for completeness (see also proposition 19 in
\cite[Sec.5.11]{Ga983},\cite[Appendix3]{Ga985}).  The equation $\V a=\V
a_0+\V n(\V a)$, with $\V a_0=-J_0^{-1}\BDpr\lis f(\V 0)$ and $\V n(\V
a)=-J_0^{-1}\BDpr_{\V a}\lis f(\V a)+J_0^{-1}\BDpr_{\V a}\lis f(\V 0)$ is
written as equation for $\V b=\V a-\V a_0$ with $\V a_0=-J_0^{-1}\BDpr_{\V
  a} \lis f(\V0)$:
\be
\eqalign{
  \V b&=%-\V a_0+
  \V n(\V a_0+\V b), \quad\ie\quad \V b=\V c+\wt{\V n}(\V b)\cr}
  \Eq{eA.1}\ee
  with $\wt{\V n}=\V n(\V a_0+\V b)-\V n(\V a_0)$
  defined in $C_{\frac12\r_0}$ if $|\V
  a_0|<\m\,\r_0$, with $\m<\frac12$, and $\V c=\V n(\V a_0)$.%-\V a_0$.
  Then the following dimensional estimates
hold:

\be\eqalign{
  |\V a_0|&\le \e_0 J_0^{-1}\le \m \r_0,\kern1.9cm {\rm if}\
  \th_0<\m \le \frac12\cr
  |\V c|\ &\le \e_0 J_0^{-1}\r_0^{-1}|\V a_0|\le\th_0\e_0 J_0^{-1}
  \equiv\th_0^2\r_0
\cr
|\wt{\V n}|\ &\le \e_0 J_0^{-1}\r_0^{-1}\,|\V b|=\th_0|\V b|,\qquad {\rm
  if}\ |\V b|<(1-\m)\r_0\cr
}\Eq{eA.2}\ee
Consider $\V b$ moving on the circle $|\V b|=\l \r_0$, $\l<1$; then:
\be
|\V b-\wt {\V n}(\V b)-\V c| \cases{\ge \l\r_0(1-\th_0-\l^{-1}\th_0^2)\cr
  \le \l\r_0(1+\th_0+\l^{-1}\th_0^2)\cr}\Eq{eA.3}
\ee
If the problem is in dimension $1$ (\ie $a,b,c,n$ are scalars) this means
that the image of the circle delimiting $C_{\l\r_0}$ is contained in the
disk with radius $\l\r_0(1+\th_0+\l^{-1}\th_0^2)$ and contains the disk
with radius $\l\r_0(1-\th_0-\l^{-1}\th_0^2)$ hence the equation has a
solution $\V b_0$ contained in the larger disk {\it if its radius is $<\r_0$
  and if the smaller radius is $>0$} (hence the latter contains the origin).

Choosing $\m=\frac12,\l=\frac13$ and assuming $\th_0<\frac1{16}$ all the
conditions are fulfilled; so the equation has a solution $\V a=\V a_0+\V
b_0=-J_0^{-1} \BDpr_{\V a} \lis f(\V a)$ in $C_{\r_0}$ (hence $|\V
a|<\th_0\r_0$) and the statement in Eq.\equ{e4.5} is proved if $\ell=1$. It
also follows that Eq.\equ{eA.1}, as an implicit equation for $\V b$ in
terms of $\V c$, gives $\V b$ analytic in terms of $\V c$.

The
multidimensional case can be proved in the same way: it might be expected
that the bound $\chi_\ell$ depends on $\ell$ but a careful examination of
the above argument not only works if $\ell>1$ (replacing the disk with a
product of $\ell$ disks) but also shows that the
constant $\chi_\ell$ can be chosen $\ell$-independent, \cite{Ga983,Ga985}.

Likewise the implicit equations discussed in Eq.\equ{e4.10},\equ{e4.11}
can be solved along the same lines, replacing disks with polydisks (see
details in propositions 20,21 in \cite{Ga983} or appendix 3 in
\cite{Ga985}.

\*

\def\SEC{Injectivity in (4.10)}
\section{\SEC}
\label{secB}
\inizA

Remark that the distance of the boundary of the polyannulus $\G_{\k_0}$ to
that of $\G_{\k_0-\d_0}$ is bounded, if $\frac12\le
e^{\pm\k_0}\le2$, below by $\frac12\d_0$ and above by
$2\d_0$. 

The injectivity follows by  integrating, along
the shortest path enclosed in $C_{\k_0-\d_0}$ connecting $\V z_1$ and $\V
z_2$,
the following difference:

\be\eqalign{
  &|\Ba'_{1,i}-\Ba'_{2,i} |=\int_{\Ba_1}^{\Ba_2}
  (\sum_j
  d\a_j \,\dpr_{\a_j} \dpr_{A'_i}\F_0(\V A',\Ba))\Big)\cr
  &\tto |\V z_1-\V z_2|\Big(1-\ell 
  \g_4\h_0\d_0^{-c_4}\Big)\cr
  &\ge \frac12|\V z_1-\V z_2|\cr}\Eq{eB.1}\ee
deduced after taking into
account the inequalities Eq.\equ{e4.8},\equ{e4.9} and determining the
constants $c_6,\g_6$.
\*

\0{\bf Acknowledgements: \it I am indebted to G. Antinucci, G.Gentile and
  I. Jauslin for comments on an earlier version. Partially supported by
  INFN, Roma1}  \*

\rhead{\thepage}\lhead{\tiny\today}

%\bibliography{0Bib}
%\bibliographystyle{acm}
%\bibliography{moser.bbl}

\begin{thebibliography}{10}

\bibitem{Wi971a}
K.~Wilson.
\newblock {Renormalization Group and Critical Phenomena. I. Renormalization
  Group and the Kadanoff Scaling Picture}.
\newblock {\em Physical Review B}, 4:3174--3183, 1971.

\bibitem{Ca966}
L.~Carleson.
\newblock On convergence and growth of partial sums of fourier series.
\newblock {\em Acta Mathematica}, 116:135--157, 1966.

\bibitem{Fe973}
C.~Fefferman.
\newblock {Pointwise Convergence of Fourier Series}.
\newblock {\em Annals of Mathematics}, 98:551--571, 197.

\bibitem{Ne966}
E.~Nelson.
\newblock A quartic interaction in two dimensions.
\newblock {\em {In {\sl Mathematical Theory of elementary particles}, ed. R.
  Goodman, I. Segal}}, pages 69--73, 1966.

\bibitem{CKN982}
L.~Caffarelli and L.~Nirenberg and L.~Kohn.
\newblock {Partial regularity of suitable weak solutions of the Navier-Stokes
  equations}.
\newblock {\em Communications on Pure and Applied Mathematics}, 35:771--831,
  1982.

\bibitem{Ga986}
G.~Gallavotti.
\newblock Quasi integrable mechanical systems.
\newblock {\em Phenom\`enes Critiques, Syst\`emes aleatories, Th\'eories de
  jauge, Proceedings, Les Houches, XLIII (1984), North Holland, Amsterdam},
  II:539--624, 1986.

\bibitem{Ga985b}
G.~Gallavotti.
\newblock Renormalization theory and ultraviolet stability for scalar fields
  via renormalization group methods.
\newblock {\em Reviews of Modern Physics}, 57:471--562, 1985.

\bibitem{WK974}
K.~Wilson and J.~Kogut.
\newblock The renormalization group and the {$\e$}-expansion.
\newblock {\em {\it The renormalization group and the {$\e$}-expansion},
  Physics Reports}, 12:75--199, 1973.

\bibitem{BG995}
G.~Benfatto and G.~Gallavotti.
\newblock {\em Renormalization Group}.
\newblock Princeton U. Press, Princeton, 1995.

\bibitem{Ga983}
G.~Gallavotti.
\newblock {\em The Elements of Mechanics (I edition);}.
\newblock Springer Verlag, New York, 1983 [I edition].

\bibitem{Ga985}
G.~Gallavotti.
\newblock Perturbation theory for classical {H}amiltonian systems.
\newblock {\em in {\it Scaling and self similarity in Physics}, Ed. J.
  Fr{\"o}hlich, Birkh{\"a}user, Boston}, pages 359--426, 1985.

\bibitem{Ko954}
A.N. Kolmogorov.
\newblock On the preservation of conditionally periodic motions.
\newblock {\em In Lecture Notes in Physics, Stochastic behavior in classical
  and quantum Hamiltonians, ed. G. Casati, J. Ford, Vol. 93, 1979}, 93, 1979.

\bibitem{Ar963b}
V.~Arnold.
\newblock {Small denominators and problems of stability of motion in classical
  and celestial mechanics}.
\newblock {\em Russian Mathematical Surveys}, 18:85--191, 1963.

\bibitem{Mo962}
J.~Moser.
\newblock On invariant curves of an area preserving mapping of the annulus.
\newblock {\em {Nachrichten Akademie Wissenshaften G{\"o}ttingen}}, 11:1--20,
  1962.

\bibitem{BGGS984}
G.~Benettin, L.~Galgani, A.~Giorgilli, and J.~Strelcyn.
\newblock {A proof of Kolmogorov's theorem on invariant tori using canonical
  transformations defined by the Lie method}.
\newblock {\em Nuovo Cimento B}, 79:201--223, 1984.

\bibitem{Ga994b}
G.~Gallavotti.
\newblock {Twistless KAM tori}.
\newblock {\em Communications in Mathematical Physics}, 164:145--156, 1994.

\bibitem{El996}
L.H. Eliasson.
\newblock Absolutely convergent series expansions for quasi periodic motions.
\newblock {\em MPEJ (Mathematical Physics Electronic Journal)}, 2, n.4:1--33,
  1986-96.

\bibitem{GG995}
G.~Gallavotti and G.~Gentile.
\newblock Majorant series convergence for twistless kam tori.
\newblock {\em Ergodic Theory and Dynamical Systems}, 15:857--869, 1995.

\bibitem{GBG004}
G.~Gallavotti, F.~Bonetto, and G.~Gentile.
\newblock {\em Aspects of the ergodic, qualitative and statistical theory of
  motion}.
\newblock Springer Verlag, Berlin, 2004.

\bibitem{GM004}
G.~Gentile and V.~Mastropietro.
\newblock {Construction of periodic solutions of the nonlinear wave equation
  under strong irrationality conditions by the Lindstedt series method}.
\newblock {\em Journal de Math\'ematiques Pures et Appliqu\'es}, 83:1019--1065,
  2004.

\bibitem{GL999}
A.~Giorgilli and U.~Locatelli.
\newblock Kolmogorov theorem and classical perturbation theory.
\newblock {\em NATO ASI series, Hamiltonian systems with three or more degrees
  of freedom}, 533:72--89, 1999.

\bibitem{Mk983}
R.~MacKay.
\newblock A renormalization approach to invariant circles in area-preserving
  maps.
\newblock {\em Physica D}, 7:283--300, 1983.

\bibitem{Ga982b}
G.~Gallavotti.
\newblock A criterion of integrability for perturbed nonresonant harmonic
  oscillators. {Wick Ordering} of the perturbations in classical mechanics and
  invariance of the frequency spectrum.
\newblock {\em Communications in Mathematical Physics}, 87:365--382, 1982.

\bibitem{CJB999}
C.~Chandre, H.~Jauslin, and G.~Benfatto.
\newblock {An Approximate KAM-Renormalization-Group Scheme for Hamiltonian
  Systems}a.
\newblock {\em Journal of statistical physics}, 94:241--251, 1999.

\bibitem{HI004}
J.~Hubbard and Y.~Ilyashenko.
\newblock A proof of {K}olmogorov's theorem.
\newblock {\em Discrete and continuous dynamical systems}, 10:367--385, 2004.

\bibitem{Ga995d}
G.~Gallavotti.
\newblock {Invariant tori: a field theoretic point of view on Eliasson's work}.
\newblock {\em Advances in Dynamical Systems and Quantum Physics, {\rm Ed. R.
  Figari, World Scientific}}, 164:117--132, 1995.

\bibitem{BGK999}
J.~Bricmont, K.~Gawedzki, and A.~Kupiainen.
\newblock Kam theorem and quantum field theory.
\newblock {\em Communications in Mathematical Physics}, 201:699--727, 1999.

\bibitem{CGJ998}
C.~Chandre, M.~Govin, and H.~R. Jauslin.
\newblock {Kolmogorov-Arnold-Moser} renormalization-group approach to the
  breakup of invariant tori in {H}amiltonian systems.
\newblock {\em Physical Review E}, 57:1536--1543, 1998.

\bibitem{Ko004}
H.~Koch.
\newblock A renormalization group fixed point associated with the breakup of
  golden invariant tori.
\newblock {\em Discrete and continuous dynamical systems}, 101:881–909, 2004.

\bibitem{Ge010}
G.~Gentile.
\newblock {Quasi-periodic motions in dynamical systems. Review of a
  renormalisation group approach}.
\newblock {\em {Journal of Mathematical Physics}}, 51:015207 (+34), 2010.

\bibitem{GG002}
G.~Gallavotti and G.~Gentile.
\newblock Hyperbolic low-dimensional invariant tori and summations of divergent
  series.
\newblock {\em Communications in Mathematical Physics}, 227:421--460, 2002.

\bibitem{GG005e}
G.~Gallavotti and G.~Gentile.
\newblock Degenerate elliptic resonances.
\newblock {\em Communications in Mathematical Physics}, 257:319--362, 2005.

\bibitem{CGGG006}
O.~Costin, G.~Gallavotti, G.~Giuliani, and G.~Gentile.
\newblock {Borel summability and Lindstedt series}.
\newblock {\em Communications in Mathematical Physics}, 269:175--193, 2006.

\end{thebibliography}
\bibliographystyle{unsrt}

\end{document}